
\input amstex
\input amssym.tex
\input amssym.def

\documentstyle{amsppt}
\magnification = \magstep1
\font\bigbf=cmbx10 scaled \magstep2

\vsize=8.55truein

\baselineskip=14pt
\overfullrule=0pt
\define\cC{{{\Cal C}}}

\define\IR{{{\Bbb R}}}

\define\IN{{{\Bbb N}}}

\define\Hom{{\text{Hom}}}
\define\out{{\text{out}}}

\vglue .5in

\centerline{\bigbf Towards a theory of classification} 

\vskip .25in   
  
\centerline{George A. Elliott}

\vskip .3truein

\midinsert\narrower\narrower

\noindent{\bf Abstract.}
The well-known difficulties
arising in a classification which
is not set-theoretically trivial---involving
what is sometimes called a non-smooth quotient---have
been overcome in a striking way in the
theory of operator algebras by the use of what
might be called a classification functor---the very 
existence of which is already a surprise. Here the
notion of such a functor is developed abstractly,
and a number of examples are considered (including
those which have arisen for various classes of operator algebras).
\endinsert

\vskip .3in
{\bf 1.}\ \
The purpose of this note is to propose
an approach to the general question of
classification.

Except in the simplest cases (sets,
vector spaces, finite simple groups!)~it is not
possible even in principle to label the
isomorphism classes of a given class of
mathematical objects in a reasonable way.
Even when this is possible, one is often
interested in more than just when two objects
are isomorphic. (For instance, one might be
interested in when one object is isomorphic to a
subobject of another.) (This is obviously important
for sets. For finite simple groups it is an
open problem.)

Perhaps, rather than labels for isomorphism
classes, what one really wants, given a
category, is a functor---distinguishing isomorphism
classes---from this category into some other, simpler,
category. (In other words, still labels for objects in
the given category, but with isomorphic objects no longer 
required to have the same

\noindent
\underbar{\hskip 2truein}

The research of the author was supported by a grant from the
Natural Sciences and Engineering Research Council of Canada. 

\vfill
\eject

\noindent
label---just isomorphic
labels!) (And maps between objects reflected by
maps---or at least formal arrows---between labels---but
reflected faintly, in the sense that certain maps will
coalesce.)

It would be natural to call such a
functor a classification functor, and the codomain
category a classifying category.

The case in which every object of the
classifying category was isomorphic to the image
of some object in the given category, and every
map between the images of two objects was
the image of a map between them---the case that
the functor was full---, would be of particular
interest. (On the other hand, the case of a
faithful functor---distinct maps in the given category taken by
the functor into distinct maps in the
codomain category---would not be of interest, the
whole point being to forget at least something!)

In this note I shall review some
examples of classification functors. Naturally, the
more concrete such a functor is, the more interesting
and useful it is likely to be. On the other hand,
one may hope that a more abstract classifying
category could also be interesting. (Recall that, in fact,
every category is concrete---it can be described as a subcategory
of the category of sets, with the maps being those  preserving 
certain operations on the sets.)

Before stating some positive results (Theorems 1 and 3 below), let
me first be a little more specific
concerning the approach of just looking at
isomorphism classes, and what the difficulties
with this approach are.

A given category, which it is desired to
classify, may or may not have other maps than
just isomorphisms. In either case, it may have
a natural topology or Borel structure---typically,
a Polish topology or a standard Borel structure.
In this case, even if one ignores maps and just
looks at objects, the quotient topology or Borel
structure on the isomorphism classes will
in general be singular. (This is one of the lessons
of the theory of operator algebras!) (The phenomenon of
non-smooth quotients was first explored by Mackey, Dixmier,
and Glimm in the setting of irreducible
representations of a C*-algebra or locally compact group.
It was later also studied by Gabriel in the setting
of indecomposable finite-dimensional representations of
directed graphs.)

In a category with only isomorphisms,
just looking at objects means forgetting
about the number of isomorphisms between two objects,
and just keeping
track of whether there is  one or not. Passing to
isomorphism classes, then, even if it destroys a given
well-behaved topology or Borel structure, may,
at least in a trivial way, be thought of as
passing by a functor to a quotient category.

If the given category includes homomorphisms
that are not isomorphisms---the most interesting
setting---the quotient category may not even
exist!

What I mean by this is that, if one identifies
arbitrary homomorphisms between two objects (instead
of just isomorphisms) whenever they differ by
an automorphism, on either the domain or codomain side, or both,
then, while this determines an equivalence
relation on the morphisms of the given category, this
is  not in general compatible with composition of
morphisms. The product of two equivalence classes
of morphisms, while it is always a union of
equivalence classes, may fail to be a single
equivalence class. 
Thus, already one may no longer
have a category (however well behaved a given
topology or Borel structure may still be). This difficulty
will persist (not to mention the possible collapse of the
topology or Borel
structure!) on passage to the isomorphism classes. (In other
words, there is no quotient category.)

For instance, the quotient category fails to
exist in this way already in the case of sets.
(The product of two non-constant maps may be
constant, whereas the product of the equivalence
classes of two non-constant maps always contains
a non-constant map.)

It is perhaps worth mentioning that this difficulty does
not arise if one restricts attention to the category
of injective maps between sets. The quotient
construction described above yields the category of
cardinal numbers, with a (unique) morphism between two
cardinals whenever the second is greater than or
equal to the first. This classifying category thus
retains the subobject information from the original
category. (Starting with the category of sets and
surjections, one also obtains a classifying category
by this construction---namely,
ignoring automorphisms. It is interesting to note
that, while the two original categories, sets with
injections and sets with surjections, would not
seem to be simply related, the classifying categories
that we have obtained for them are related in a very
simple way: one is just the dual of the other---i.e.,
the category with the same objects and all the arrows
reversed.) (This is not the case  for the original
categories---as can be seen by just counting numbers of arrows, for 
instance when the domain set is a single point.) 
(It is fortunate that the constructions just
described work, as, as is well known, cardinal
numbers cannot just be defined as equivalence classes
of sets---also the order relation is needed.)

The case of vector spaces is very similar to the case
of sets (and may be essentially reduced
to that case by choosing bases).

It turns out that a somewhat similar
difficulty occurs with the category of finite simple
groups (even with injective maps). Since there is an automorphism of the group
$A_6$ taking the permutation (123) into the permutation
(123)(456) (I am indebted to J.~B.~Olsson for this
calculation), this does not extend to an automorphism
of $A_7$ (the automorphism group of which, in contrast
to the case of 
$A_6$, is just $S_7$). The product of the equivalence classes
(modulo automorphisms) of the
canonical embeddings of $A_3$ in $A_6$ (as the permutations
of the first three symbols) and of $A_6$ in $A_7$ (as the
permutations of the first six symbols), because of this
automorphism of $A_6$, is therefore
strictly larger than the equivalence class of the product,
the canonical embedding of $A_3$ in $A_7$. (It also contains
the equivalence class of the embedding of $A_3$ in $A_7$ 
with multiplicity two, i.e.,
taking the permutation (123) into (123)(456).)

Somewhat as for sets (or for vector spaces), the situation
for finite simple groups---indeed, for arbitrary groups---can
be salvaged---in this case, by
throwing away (dividing out by) fewer than all
automorphisms, namely just the inner ones. (This does
not amount to anything in the commutative case; the
theory proposed is a purely noncommutative one.)

In other words, if one identifies two group homomorphisms
if they differ by an inner automorphism---on either the
domain side or the codomain side or both---then, not
only does one obtain an equivalence relation, but also,
it is compatible with composition: the
product of two equivalence classes is again an
equivalence class. (This is because the composition
of a morphism with an inner automorphism on the
domain side is equal to the composition of the
same morphism with another inner automorphism on
the codomain side.)

Clearly, the resulting functor is a classification
functor---and furthermore, a similar
construction works in other settings, for instance
in the category of rings (where by an inner automorphism
of a ring is meant one determined by an invertible element 
of the ring obtained by adjoining a unit).

Let us formalize this construction.
\medskip

{\bf Theorem.}\
{\it Let $\cC$ be a category with a notion
of inner automorphism, satisfying the axiom that
the composition of an arbitrary morphism with (what
we shall call) an
inner automorphism, on the domain side, is equal
to the composition of the same morphism with 
another inner automorphism on the codomain side (just
as recalled above for groups and for rings). 
Note that, given sets of inner automorphisms, in this
sense, simultaneously for all objects in the category,
the subgroups generated by these sets of automorphisms also
satisfy the axiom, and in particular 
are normal subgroups.
Overall, these subgroups (which we shall refer
to as the inner automorphism groups)
form what we might
refer to as a compatible family of normal subgroups
of the automorphism groups.

It follows that the category $\cC^{\text{out}}$, the objects of
which are the same as those of $\cC$, and the morphisms of
which are those of $\cC$ considered modulo
inner automorphisms, is a classifying category
for $\cC$.}
\medskip

{\it{Proof.}}\
The main point is that (cf.~above) $\cC^{\text{out}}$
is a category. It is immediate that,
if a map in $\cC$ is
invertible in $\cC^{\text{out}}$ then it is invertible.
(And so the canonical
functor from $\cC$ to $\cC^{\text{out}}$ is a classifying
functor.)
\medskip

{\bf 2.}\ \
A particular case of Theorem 1 is the
category of (non-zero) finite direct sums of
matrix algebras over the complex numbers,
considered as  *-algebras---i.e., with 
*-homomorphisms as maps---with the operation of taking
the adjoint, or conjugate transpose,
of a matrix as the *-operation. In this case
(and for that matter in the category of all 
C*-algebras, of which this is a subcategory),
one has the compatible family of normal
subgroups of the automorphism groups consisting
of the inner automorphisms, i.e., the automorphisms
determined by unitary elements of the *-algebra obtained
by adjoining a unit.

As pointed out by Bratteli in [1], in this case  the
classifying category constructed in Theorem 1
has a very simple description, combinatorial in
nature. The objects (non-zero finite direct sums
of matrix algebras) may be viewed as (i.e., labelled
precisely by) the finite column vectors, of
arbitrary (non-zero) length, the coordinates
of which are strictly positive integers. The
morphisms between two objects, or vectors, may then be
viewed as the rectangular matrices with
positive (not necessarily strictly positive) integers
as entries, multiplying the first vector into
either the second vector (if the map is unital),
or a vector with smaller coordinates. (Here the multiplication by
the matrix is understood to be on the left, and
the numbers of columns and rows of the
rectangular matrix must therefore be equal
respectively to the numbers of coordinates of
the domain and codomain vectors.) In this
description, according to the computation of Bratteli,
composition of  *-algebra morphisms---modulo
inner automorphisms---corresponds to 
multiplication of rectangular matrices.

In slightly different words, if the set of
column vectors and rectangular matrices described
above is considered with its natural structure
as a category (described above), then one obtains
an exact replica  (up to equivalence of categories)
of the classifying category
of Theorem 1, in the case of the category of 
*-algebras under consideration. What this
comes down to is that if one considers two
single full matrix algebras, then there is at most
one unital morphism from the first to the second, up to
unitary equivalence---and this exists exactly when
the order of the second matrix algebra is an
integral multiple of the order of the first one
(sometimes called the multiplicity of the
embedding). 
(Similarly, a non-unital morphism is also determined
up to unitary equivalence by its multiplicity---defined
by cutting down by the image of the
unit and so reducing to the unital case. 
The multiplicity can be any positive integer the product of
which with the order of the first matrix algebra is 
less than or equal to the order of the second matrix algebra.)

The Bratteli matrix, in the case
of a  general pair of algebras in the category
under consideration, 
 just keeps track of the
multiplicities of what might be called the
partial maps, from the individual minimal
direct summands of the domain algebra to
those of the codomain algebra.
\medskip

{\bf 3.}\ \
In fact, Bratteli was interested in a
larger category, 
namely, in the category of all C*-algebra
inductive limits of sequences in the category of
C*-algebras considered above (non-zero finite
direct sums of matrix algebras over the complex numbers).
(Equivalently, Bratteli considered the category of 
C*-algebras obtained as the closure
of an increasing sequence of sub-C*-algebras belonging to the
category of Section 2.)

While Theorem 1 is applicable to this category
also, in fact the classifying category arising
in this way suffers from one of the defects described
in Section 1: it is singular. (The
group of inner automorphisms of a general C*-algebra in this
category is not a closed subgroup of the
group of all automorphisms in its natural topology, and
so the quotient is not Hausdorff.) (This remark applies
also to many other categories, for instance, infinite groups.)

What Bratteli did instead, circumventing this
difficulty, was, given a sequence of C*-algebras in the category
of Section 2, to look at the sequence in
the classifying category of this category given by
Theorem 1. In his picture, this was a
diagram consisting of a whole sequence
of column vectors, each one connected to the next by a
rectangular matrix, as
described in Section 2. This is now called
a Bratteli diagram.

What Bratteli observed,
to a certain extent implicitly, was that
in a natural way the Bratteli diagrams form a category, and that
if, for every C*-algebra in his category (which he
called the approximately finite-dimensional, or  AF,
C*-algebras), one just chooses a Bratteli diagram
(from a particular representation of this algebra as
an inductive limit), then one obtains a classification
functor. (In fact, Bratteli considered only isomorphisms,
but his considerations can be extended in a natural way
to embrace
arbitrary morphisms.)

What I propose to do here is to take
Theorem 1 seriously for the larger category,
and indeed also for
many other categories---e.g., all separable
C*-algebras, and all countable groups.

It turns out that it is possible to desingularize
Theorem 1.
\medskip

{\bf Theorem.}\
{\it{Let $\cC$ be a category with a notion
of inner automorphism, i.e., a compatible family
of normal subgroups of the automorphism groups
as described in Theorem 1. Suppose that for each
pair of objects the set of morphisms between these
objects is endowed with
a  complete metric space structure, and that
the following two compatibility properties with
regard to composition of morphisms hold.

First, for any three objects, composition
of morphisms from the first object to the second with
morphisms from the second object to the third is 
a (jointly) continuous map into the space of morphisms
from the first object to the third. (This property
pertains to the topology, not the metric itself.)

Second,  for any two objects, and for
a fixed inner automorphism of the second object,
composition with this (on the codomain side) is an
isometry  from the space of all morphisms from
the first object to the second  onto itself.

If follows from the first axiom (continuity) that
the quotient structure $\cC^{\text{out}}$, the objects of which are
the same as those of $\cC$ (and as those of the category
$\cC^{\text{out}}$ of Theorem 1), and the morphisms of which are
the closures of the equivalence classes of morphisms
of $\cC$ modulo inner automorphisms (alternatively---as by
continuity the closure of an equivalence class is a union 
of equivalence classes---the
closures of the morphisms of the category $\cC^{\text{out}}$---in
the quotient topology in which
points are not necessarily closed), and for which
the product of two morphisms is defined as the
closure of the product of the corresponding two
closed sets of morphisms of $\cC$ (alternatively, as
the closure of the product of the corresponding
closures of single morphisms in $\cC^{\text{out}}$)---by
continuity this is the
closure of a single equivalence class
of morphisms of $\cC$ (i.e., a single morphism of $\cC^{\text{out}}$,
namely, the product of two morphisms generating the two
point closures in question), and therefore it is a morphism---is
a category. Furthermore, the quotient map is a functor.

It follows from the second axiom (or, rather, the two axioms together)
that the natural functor from the category $\cC$ to the
category $\overline{\cC^{\text{out}}}$
(i.e., the quotient map) is a classification functor.
(In other words, it distinguishes isomorphism classes.) It is in
fact what might be called a strong classification functor,
in the sense that isomorphisms lift to $\cC$ from the classifying 
category $\overline{\cC^{\text{out}}}$.}}
\medskip

{\it  Proof.}\
The main point is still that $\overline{\cC^{\text{out}}}$ is
a category. (The proof that the natural functor from
$\cC$ to $\overline{\cC^{\text{out}}}$
distinguishes isomorphism classes is,
as will be seen, not new.) 

It must be checked that
composition of morphisms (as defined in the
statement of the theorem) is associative. (Strictly
speaking, this must also be checked for $\cC^{\text{out}}$: the
product of the equivalence classes of three
morphisms in $\cC$, in a fixed order, of course, but
grouped in either way, is just the equivalence
class of the product of the three morphisms; this
may be seen immediately by, roughly speaking, just
moving all inner automorphisms through to the
codomain side.)

Once it is noted that composition of morphisms
in $\cC^{\text{out}}$---i.e., equivalence classes of morphisms in
$\cC$---is associative, it follows immediately by continuity
of multiplication that
composition of morphisms
in $\overline{\cC^{\out}}$---i.e., closures of equivalence classes
in $\cC$---is associative: just as one sees (by continuity) that
the closure of the product of the closures of
two equivalence classes is just the closure of the
product of the equivalence classes themselves (and in particular
is the closure of a single equivalence class), so also
one sees that the two sets
involved in the law of associativity for closures
of equivalence classes---with multiplication of two
such closures the closure
of the product---equivalently, the closure of the
product of the equivalence classes themselves---are
equal (each one equal to the closure of the product of
all three equivalence classes in question---of course, this
uses associativity of the product of equivalence classes).

It is also immediate, starting from the
continuity of multiplication in $\cC$ and the
functoriality of the quotient map from $\cC$ to $\cC^{\out}$,
that the quotient map from $\cC$  to $\overline{\cC^{\out}}$ is a
functor. Indeed, functoriality from $\cC$ to 
$\cC^{\out}$ just
says that the equivalence
class of the product of two arrows in $\cC$ is
the product of the equivalence classes, and by
continuity of multiplication this implies that the
closure of the equivalence class of the product of two
arrows is the closure of the product of the
closures of the equivalence classes, which is the
desired functoriality.

It is interesting to note that, so far, besides
continuity of composition of morphisms (i.e., arrows) joining a
fixed triple of objects (from the first to the second
and the second to the third)---and of course associativity of
this composition---the only thing that has been
used, to obtain that $\cC^{\out}$ and (hence) 
$\overline{\cC^{\out}}$ are categories and that the
natural maps $\cC \to \cC^{\out}$ 
and (hence) $\cC\to \overline{\cC^{\out}}$ are functors,
is that the product of two equivalence classes
is again an equivalence class.
Whereas to prove that $\cC\to \cC^{\out}$
is a classification functor, also nothing more
is needed, to prove that $\cC\to\overline{\cC^{\out}}$ is a
classification
functor seems to require the full force of the hypotheses,
i.e., that the equivalence classes
derive from neglecting the so-called inner automorphisms
(which of course in particular means that the product
of two equivalence classes is again an equivalence
class), and that the topology on the set of
morphisms from each fixed object to another
one derives from a metric, which is assumed to be both complete
and invariant under composition with a fixed
inner automorphism of the codomain object.

It may also be of some interest to note that
the stronger hypotheses, crucial for the
second statement of the theorem, also have two
incidental consequences which one might think related
to the first statement, that one has a
quotient category and a functor to it, but do not
appear to be so related: First, the quotient map
$\cC\to \cC^{\out}$
is open (when each set of morphisms between a pair of
objects in $\cC^{\out}$ is given the quotient topology), and
indeed the quotient map $\cC\to {\overline{\cC^{\out}}}$ is 
open (again with respect
to the quotient topology---equivalently, 
the quotient map $\cC^{\out} \to {\overline{\cC^{\out}}}$ is 
also open---in fact the               
saturated open sets of morphisms in $\cC$ 
with respect to the map $\cC\to \cC^{\out}$ are already 
saturated with
respect to the map $\cC\to \overline{\cC^{\out}}$! Second, the
closures of distinct equivalence classes of morphisms
in $\cC$ are either equal, or disjoint. (Of course,
even without the stronger hypotheses, the most that can
happen if the closures of two equivalence classes are
neither equal nor disjoint is that one is contained
in the other, properly, but this can 
presumably happen. Interestingly, this does not create
any difficulty in the definition of the category structure of
$\overline{\cC^{\out}}$.)

The fact that the functor in question, that to
each morphism of $\cC$ associates the morphism in
$\overline{\cC^{\out}}$ consisting of the closure of the equivalence
class of this morphism in $\cC$ modulo inner
automorphisms (these automorphisms defined axiomatically in the
statement of the theorem), is a classification functor
(i.e., distinguishes isomorphism classes), depends on
a sequential approximate intertwining argument which
was developed first in [1] and [8] in the
case of exact intertwinings, and in [9] in the
case of approximate intertwinings (i.e., approximately
commutative diagrams intertwining two
sequences). This technique has been used many times
since, indeed, on virtually every occasion that an
isomorphism theorem for C*-algebras has been established. 
(It might almost be omitted, so many times has it been used!
The basic ingredients were already described abstractly in
[9]---see Theorems
2.1 and 2.2 of the that article.
Surveys of the C*-algebra isomorphism results are given
in [12] and in [29].) The main purpose of the
present article is to suggest that it might
be of serious interest beyond the setting of
C*-algebras. (For instance, for von Neumann algebras!)
(And also, for countable (non-abelian) groups.)
(Some observations in this direction are reported below, in
Sections 4 and 5.)

Let $a$ and $b$ be objects in $\cC$, and
suppose that they are isomorphic in $\overline{\cC^{\out}}$. 
Let $f$ be an isomorphism between $a$ and $b$ in the
category $\overline{\cC^{\out}}$, and let us show that $f$ is
the image of an isomorphism in $\cC$. (In order to
prove that the functor $\cC\to \overline{\cC^{\out}}$ is a
classification functor, it
appears to be expedient to prove that it is
a strong classification functor. Indeed, only in
[28] (and later, in a similar way, in [15]) 
has a classification functor been obtained
without showing that it is a strong classification
functor. Recently, in [7], the
functor shown to be a
classification functor in [28] was shown in fact
to be a strong classification functor. In this
connection, the proof of Theorem 1 shows that
the functor $\cC\to\cC^{\out}$ is not only a strong
classification functor, in the sense that every
isomorphism in the codomain (classifying) category is the
image of an isomorphism in the domain category,
but is what might be called super-strong, in
the sense that any morphism in the domain
category mapping into an isomorphism in
the classifying category must in fact already be
an isomorphism!)

Choose a morphism $f_1$ in $\cC$ mapping
to $f$ by the functor, and a morphism
$g_1$ in $\cC$ mapping to the inverse of $f$.
Consider the (non-commutative!) diagram
$$\eqalign{
&  a           \to  a  \to  a  \to  \cdots  \to  a \cr
& \hskip -.05truein   \downarrow  \ \nearrow \  \downarrow \   \nearrow \  \downarrow \ 
\nearrow\cr
& b \to \ b \to  b  \to   \cdots \to b\cr}
$$
in which all the horizontal arrows are the
identity map, respectively for $a$ and for $b$, and
the downwards and upwards arrows are $f_1$
and $g_1$ respectively (each repeated 
infinitely often). Let us modify this diagram,
by changing each of the downwards and upwards
arrows by inner automorphisms, to make it
approximately commutative---in the natural sense
described in a special case in Section 2 of [9]
and, in the present abstract setting, implicitly 
below. For convenience of notation, let us relabel the
downwards maps, $f_1, f_1, \cdots$, as $f_1, f_2, \cdots$,
anticipating that they will be changed.
Similarly, let us relabel the
upwards
maps $g_1, g_1, \cdots$
as $g_1, g_2, \cdots$ (of course, to begin with,
all the same).

Consider first the upper left hand  triangle
in the diagram. It provides two routes from  $a$
to $a$, which we might refer to as ``across" and
``down-up", which agree exactly
in $\overline{\cC^{\out}}$, by hypothesis, and in other words are
approximately equal in $\cC$ modulo inner automorphisms
(to within an arbitrarily close degree of approximation).
Since multiplying by an inner automorphism 
(on the codomain side)
preserves distances, in the space of morphisms from $a$ to
$a$, it is possible
to multiply just one of the morphisms by an inner automorphism 
to get one arbitrarily close to the other one.
Doing this to the map ``down-up", i.e., to
$f_1 g_1$ (we are using category theory notation for
composition of arrows, with the first on the left),
we obtain that $f_1g_1 h$ is at distance at most
$2^{-1}$ from ``across", i.e. from $\text{id}_a$, for  some
inner automorphism $h$ of $a$. 
Replacing $g_1$ by $g_1 h$,
but keeping the same notation, we then  have 
$$d(f_1 g_1, \text{id}_a) \le 2^{-1},$$
where $d$ denotes the invariant metric.

Similarly, considering the second triangle (the
lower left hand one), with  (the new) $g_1$ as the 
map ``up", and $f_2$ as the map ``down", 
noting that the two
routes ``across" and ``up-down" are (still) exactly equal
in $\overline{\cC^{\out}}$, and therefore approximately equal
in $\cC$ modulo inner automorphisms, and using that
the metric on the space of morphisms from $b$ to $b$ in
$\cC$ is invariant under multiplying
(on the codomain side) 
by inner
automorphisms, we obtain that $g_1 f_2 k$,
the map ``up-down", is within
distance $\epsilon_2$ of $\text{id}_b$, the map ``across",
for some inner automorphism $k$ of $b$, where $\epsilon_2$ is to
be specified.  Replacing $f_2$ by
$f_2 k$ (but keeping the same notation), we now have
``up-down" close to ``across" in the second triangle:
$$d(g_1 f_2, \text{id}_b)\le \epsilon_2.$$

Continuing in  this way, changing the right hand
non-horizontal arrow in each triangle in turn
(but not the left hand one, which was the
right hand one of the previous triangle---and therefore
not interfering with the approximate commutativity
of any previous triangle), we arrive at a new
choice of the non-horizonal arrows in the diagram,
$f_1, f_2, \cdots$ and $g_1, g_2, \cdots$, agreeing with
the original choice up to inner
automorphisms, and such that, for each $n=1, 2, \cdots,$
$$
\hskip .2truein  d(f_n g_n, \text{id}_a)  \le \epsilon_{2n-1}$$
and
$$ d(g_n f_{n+1}, \text{id}_b) \le \epsilon_{2n},
$$
where $\epsilon_1 =2^{-1}$, and $\epsilon_{2n-1}$ and
$\epsilon_{2n}$ are to be specified.

We wish to choose $\epsilon_2, \epsilon_3, \cdots$ in
such a way that the sequences of morphisms $(f_n)$
and $(g_n)$ converge, say to $f_\infty$ and $g_\infty$, 
from $a$ to $b$ and from $b$ to $a$, respectively.
Of course, this will probably require that
$\epsilon_k$ tend to zero, but if we actually ensure that this
holds, then, in addition, necessarily $f_n g_n$ converges to
the identity for $a$, $\text{id}_a$, and $g_n f_{n+1}$ converges to
$\text{id}_b$. By continuity of multiplication,
this implies that
$$f_\infty g_\infty =\text{id}_a,\ g_\infty f_\infty =\text{id}_b.$$

Since $f_\infty$ and $g_\infty$ still map into $f$ and $g$ in
$\overline{\cC^{\out}}$ (being limits of elements in the
equivalence classes of
$f_1$ and $g_1$ respectively in $\cC$), this shows that
the given isomorphism $f$ from $a$ to $b$
in $\overline{\cC^{\out}}$ lifts to an isomorphism in $\cC$, as
desired.

Consider the following choice of the sequence
$\epsilon_1, \epsilon_2, \cdots$. Keep $\epsilon_1=2^{-1}$;
we are embarking on what
is famously known as ``a $2^{-n}$ argument". (``The
construction of a summably Cauchy sequence" 
might be a clearer description!) Choose $\epsilon_2$, 
using continuity of multiplication, small enough that
$$d(f_1 (g_1 f_2), f_1)\le 2^{-2}.$$
(Recall that $d(g_1 f_2, \text{id}_b)\le \epsilon_2$, and
that $f_1$ is fixed; in fact also $g_1$ is fixed, and
only $f_2$ has
to be chosen suitably to ensure
that $\epsilon_2$ is small---and
this is crucial in obtaining the original sequence of estimates,
one for each $\epsilon_n$---but this is now no longer needed.)
Similarly, choose $\epsilon_3$ small enough that
$$d(g_1(f_2 g_2), g_1)\le 2^{-3}.$$
(Recall that $d(f_2 g_2, {\text{id}}_a)\le \epsilon_3$, and that, before
we came to consider
the final choice of $f_2$ and $g_2$, we had already fixed on a
choice of $g_1$; in fact, the inequality
$d(f_2 g_2, {\text{id}}_a)\le \epsilon_3$
was negotiated without changing an
earlier---final!---choice of $f_2$, and to
ensure that $\epsilon_3$ is small it is necessary only
to modify $g_2$---and, again, this is crucial in obtaining the
estimates involving $\epsilon_1, \epsilon_2, \cdots$ in the
first place---which depended on
the choice at each stage not affecting earlier 
estimates---but it is not
needed now.)
\medskip

Continuing step by step in this way, we obtain a sequence
$(\epsilon_1, \epsilon_2, \cdots)$ such that for each
$n=1, 2, \cdots,$
$$d(f_n(g_n f_{n+1}), f_n)\le 2^{-(2n-1)},$$
and
$$d(g_n (f_{n+1} g_{n+1}), g_n)\le 2^{-2n}.$$

Next, revisiting the choice of the sequence
$(\epsilon_1, \epsilon_2, \cdots)$, let us revise it slightly,
making each $\epsilon_k$ in turn smaller, sufficiently small that,
from
$d(f_n g_n, \text{id}_a)\le \epsilon_{2n-1}$
and $d(g_n f_{n+1}, \text{id}_b)\le\epsilon_{2n}$,
it follows by continuity that, for each $n=1, 2, \cdots$, 
$$d((f_n g_n)f_{n+1}, f_{n+1})\le 2^{-(2n-1)},$$
and
$$d((g_n f_{n+1}) g_{n+1}, g_{n+1})\le 2^{-2n}.$$

\noindent
This requires some comment, since, after the choice of
$\epsilon_{2n-1}$, the morphism $f_{n+1}$ will be
changed, in order to ensure that
$d(g_n f_{n+1},\ \text{id}_b) \le \epsilon_{2n}$!
But this (single) change consists only in multiplying 
by an inner automorphism, which will not affect
the first inequality above.
Similarly, the modification of $g_{n+1}$ by an
inner automorphism after the choice of $\epsilon_{2n}$,
in order to ensure that $d(f_{n+1} g_{n+1}, \text{id}_a)\le
\epsilon_{2n+1}$, does not
affect the second inequality above.

Matching the appropriate pairs of inequalities, we
obtain
$$d(f_{n+1}, f_n)\le 2(2^{-(2n-1)})=2^{-2n+2},$$
and
$$d(g_{n+1}, g_n)\le 2 (2^{-2n}) =2^{-2n+1}.$$

In particular, the sequences $(f_n)$ and $(g_n)$ are
(summably!) Cauchy, as desired. As shown
above the limits
are necessarily the inverses of each other, and
give rise to the map $f$ and its inverse in
$\overline{\cC^{\out}}$.
\medskip

{\bf 4.}\ \
{\bf Examples}
\medskip

{\bf 4.1.}\ \
{\bf Countable groups.}\ \ Consider the
category of countable (discrete)
groups, and for each object consider the normal
subgroup of the automorphism group consisting of the
inner automorphisms in the usual sense. For each
object $G$ choose a
numbering of its elements, i.e., a bijection 
$G\ni g\mapsto n(g)$ of $G$ onto either $\IN$
or an initial segment of $\IN$, and note that for 
each object $H$ the formula (involving the
Kronecker delta)
$$d(\varphi, \psi)=\sum_{g\in G} 2^{-n(g)}\delta_{\varphi(g), \psi(g)}$$
defines a metric on the set $\Hom (G, H)$ of group
homomorphisms from $G$ to $H$. (For each $g\in G$ the
quantity $d_g (\varphi, \psi)=\delta_{\varphi (g), \psi(g)}$
is already a pseudo-metric on
$\Hom (G, H)$, i.e., satisfies the triangle inequality.)

The resulting family of normal subgroups is a 
compatible one in the sense of Theorem 3
(and Theorem 1), and the
resulting family of metrics on the sets of morphisms
between pairs of objects satisfies the two axioms 
of Theorem 3.
(The underlying topology is just the topology of
pointwise convergence in the discrete topology;
composition of morphisms, between two fixed pairs of objects,
is easily seen to be continuous in this topology.
If $\varphi, \psi\in \Hom (G, H)$ and 
$\rho$ is an automorphism of $H$---not necessarily
inner---then
$\delta_{\rho\circ \varphi (g), \rho\circ\psi (g)}=
\delta_{\varphi(g), \psi(g)}$
and so
$$\eqalign{
d(\rho\circ\varphi, \rho\circ\psi) & =
\sum_{g\in G} 2^{-n(g)} \delta_{\rho\circ\varphi(g), \rho\circ\psi(g)}
\cr               
& = \sum_{g\in G} 2^{-n(g)}\delta_{\varphi(g), \psi(g)}\cr
& = d(\varphi, \psi);\cr}
$$
in other words,  composition with $\rho$ is an isometry.)
\medskip

{\bf 4.2.}\ \
{\bf Countably generated algebras.}\ \ Consider the
category of countably generated (not necessarily unital) 
algebras over a fixed field,
and for each object consider the normal
subgroup of the automorphism groups consisting
of the inner automrophisms in the usual sense
(i.e., those automorphisms the unique extension of
which to the algebra with unit adjoined is
determined by an invertible element). For each
object $A$ choose a generating sequence 
$(a_n)_{n\in \IN}$, and note that for each
object $B$ the formula
$$d(\varphi, \psi)=\sum_{n\in\IN} 2^{-n}\delta_{\varphi(a_n),
\psi(a_n)}$$
defines a metric on the set $\Hom (A, B)$ of
algebra homomorphisms from $A$ to $B$. (For each
$a\in A$, the quantity $d_a (\varphi, \psi)=
\delta_{\varphi(a), \psi(a)}$ is
already a pseudo-metric on $\Hom (A, B)$, i.e.,
satisfies the triangle inequality.)

The resulting family of normal subgroups
is a compatible one in the sense of Theorem 3,
and the resulting family of metrics on the sets
of morphisms between pairs of objects satisfies the
two axioms of Theorem 3. (As before, the underlying
topology is pointwise convergence (on the domain object)
in the discrete topology (of the codomain object), 
and it follows immediately that composition
is continuous. The isometry property holds since, as
before, $\delta_{\rho\circ\psi(a), \rho\circ\psi(a)} =   
\delta_{\varphi(a), \psi(a)}$ for any
$\varphi, \psi\in \Hom (A, B)$, any automorphism $\rho$---inner
or not---of $B$, 
and any
$a\in A$.)
\medskip

{\bf 4.3.}\ \
{\bf Separable C*-algebras.}\ \ Consider the
category of separable C*-algebras (not necessarily unital),
and for each object consider the normal subgroup of the
automorphism group consisting of the inner automorphisms
in the usual sense (i.e., those automorphisms
determined by a unitary element of the C*-algebra
obtained by adjunction of a unit). For each
object $A$ choose a generating sequence $(a_n)_{n\in\IN}$
of elements of norm at most one, and note that
for each object $B$ the formula
$$d(\varphi, \psi)=\sum_{n\in \IN} 2^{-n}\|\varphi(a_n)-
\psi(a_n)\|$$
defines a metric on the set $\Hom (A, B)$ of
morphisms from $A$ to $B$ (in the
usual sense of *-homomorphisms---recall that
*-homomorphisms between C*-algebras are
norm-decreasing, and so the sum is finite).

The resulting family of normal subgroups is a
compatible one in the sense of Theorem 3, and the
resulting family of metrics on the sets of morphisms
between pairs of objects satisfies the two axioms
of Theorem 3. (The underlying topology is again
the topology of pointwise convergence, now with
respect to the norm topology on the codomain
object; composition of morphisms is again
continuous, as follows from the facts that
pointwise convergence implies uniform convergence on
compact subsets. If $\varphi, \psi \in \Hom (A, B)$ and
$\rho$ is an automorphism of $B$---inner or not---then
$\|\rho\circ \varphi(a) -
\rho\circ \psi (a)\| = \|\varphi(a)-\psi(a)\|$
for each $a\in A$ (and in particular for $a_n$ for
each $n$) and as before it follows that composition with 
$\rho$ is an isometry.)

(A natural generalization of this example is the category  
of separable C*-algebras together with actions of a fixed
locally
compact group---a single object consisting of a C*-algebra
together with an action of the group---with as morphisms
those C*-algebra morphisms which are compatible with the
actions, and with as inner automorphisms
those which are inner in the sense
defined above---as C*-algebra automorphisms.
As metric on the space of morphisms between two objects,
one may just take the metric relative to a choice of
a dense sequence in the unit ball of each separable
C*-algebra as defined above. It would be interesting
to consider the even more general setting in which the  
morphisms are only required to be compatible
with the group actions up to a cocycle---isomorphism
of two actions in this category would then be what is 
known as cocycle conjugacy; work of Evans and
Kishimoto and also recent work of Katsura and Matsui
involves an intertwining argument in this setting.) 

\medskip

{\bf 4.4.}\ \
{\bf Countably generated Hilbert C*-modules with
embeddings.}\ \ Consider the
category of countably generated right Hilbert C*-modules
over a fixed C*-algebra $A$ (see e.g.~[20]), with
$A$-valued inner product preserving $A$-module maps
(not necessary adjointable) as
morphisms, and for each object consider the normal subgroup of the
automorphism group consisting of what might be
called the inner automorphisms (i.e., those
automorphisms arising from unitary elements of
the C*-algebra of compact endomorphisms with 
unit adjointed). For each object $X$ choose a generating
sequence $(\xi_n)_{n\in \IN}$ of elements of norm at most one,
and note that for each object $Y$ the formula
$$d(\varphi, \psi)=\sum_{n\in \IN} 2^{-n} \| \varphi(\xi_n)-
\psi (\xi_n)\|$$
defines a metric on the set $\Hom (X, Y)$ of morphisms
from $X$ to $Y$ (as defined above---note in particular
that morphisms are isometric and so the sum is finite).

The resulting family of normal subgroups is a
compatible one in the sense of Theorem 3, and the
resulting family of metrics on the sets of morphisms
between pairs of objects satisfies the two axioms of
Theorem 3. (The compatibility follows from the fact that
compact endomorphisms of a closed submodule of
a Hilbert C*-module extend canonically to compact
endomorphisms of the whole module.
As in the preceding example, the topology
underlying the metric is that of pointwise
convergence in norm, and continuity of
composition follows in the same way. The isometry
property
is also seen in the same way: if $\varphi, \psi \in
\Hom (X, Y)$ 
are morphisms (embeddings) 
as above, and $\rho$ is an automorphism of $Y$---inner or
not---, then 
$\|\rho\circ\varphi(\xi)-\rho\circ\psi(\xi)\|=
\|\varphi(\xi)-\psi(\xi)\|$ for each
$\xi\in X$ (and in particular for $\xi_n$ for each $n$), 
and it follows that composition with $\rho$ is an isometry.)
\medskip

{\bf 4.5.}\ \
{\bf von Neumann algebras with separable predual.}\ \ Consider
the category of von Neumann algebras with separable
predual, with a specified choice of faithful normal state for
each von Neumann algebra (this might
be called the category of pointed von Neumann algebras),
and with morphisms the (unital) normal *-homomorphisms from one
von Neumann algebra to another, 
taking the centre of the first into  (a subalgebra of)
the centre of the second
(automatic in the factor case),
which are compatible with the chosen pair
of states in the strongest sense---i.e., also intertwining
the modular automorphism groups.
(In particular just simple  intertwining of the states
implies that 
the homomorphism is 
injective.) For each object, consider the normal
subgroup of the automorphism group---i.e., the group
of *-automorphisms of the von Neumann algebra
commuting with the modular automorphism group of the
specified faithful normal state---consisting of those
automorphisms which are
inner in the usual sense---they are then determined
by unitaries fixed by  the modular automorphisms
up to central multiples.
For each object $M$---let us not
mention the given faithful normal state explicitly, but denote
the corresponding pre-Hilbert space norm by $\| \cdot\|_2$---choose
a generating sequence $(x_n)_{n\in \IN}$ of elements of
(operator) norm at most one, and note that for each
object $N$ the formula
$$d(\varphi, \psi)=\sum_{n\in \IN} 2^{-n}\| \varphi(x_n)-
\psi(x_n)\|_2$$
defines a metric on the set $\Hom(M, N)$ of morphisms
from $M$ to $N$ (recall that by assumption morphisms
preserve the canonical norm $\| \cdot\|_2$ and hence, since
$\| x_n\|_2 \le \|x_n\| \le 1$, the sum is finite).

The resulting family of normal subgroups is
a compatible one in the sense of Theorem 3,
and the resulting family of metrics on the sets
of morphisms between pairs of objects satisfies
the two axioms of Theorem 3.
(The underlying topology is again the topology of pointwise
convergence, now with respect to the strong operator
topology on the codomain object---which, on
sets bounded in the operator norm, is what the
norm $\| \cdot \|_2$ gives rise to. However, this
description is not sufficient to prove continuity of
multiplication, as continuity of composition of
arbitrary *-homomorphisms in this topology
presumably does not hold; it is necessary to
use the invariance of the norm $\| \cdot \|_2$ under morphisms
as at present defined. In this setting the proof is
very much the same as before---to prove continuity
of multiplication in the topology of pointwise convergence
with respect to the norm $\| \cdot \|_2$, it is enough to note
that the topology is sequentially determined,
and that pointwise convergence of morphisms, 
with respect to the norm $\| \cdot \|_2$, 
implies uniform convergence on each subset
totally bounded in the norm $\| \cdot \|_2$ 
(and in particular on each convergent sequence).)

The subcategory of finite
factors, with the specified state the trace, is particularly
simple as the morphisms are arbitrary (unital) normal
*-homomorphisms, and the specified normal subgroup is
the group of all inner automorphisms.
\medskip

{\bf 5.}\ \
{\bf Concrete description of the abstract classifying
category}

In certain cases the abstract classifying
category of Theorem 3 has a relatively simple concrete
form---and sometimes the classification functor when
translated into these terms can be recognized.
\smallskip

{\bf 5.1.}\ \
A good example from this point of view is the
class of countable groups obtained as inductive
limits of sequences of finite products of alternating
groups---on five or more symbols so that they are simple, and so
the maps
between building blocks (the finite products)
are determined up to inner automorphisms---or rather,
up to automorphisms determined by permutations, not necessarily even, on 
each simple component---just by
the multiplicities of the partial maps from the
various components of the domain building block to the
various components of the codomain building block. 
(Note that, as is easily seen by induction on the number of symbols of the
domain group, any homomorphism
from an alternating group on 
five or more symbols into a larger 
alternating group is determined up to a permutation,
not necessarily even,
by an integer greater than or equal to zero which might
be called the multiplicity.)

In order to ensure that the automorphism relating two
maps with the same multiplicities (for all partial maps) may be chosen
to be inner, i.e., to arise from an even permutation
on each alternating group component of the codomain group, 
we must
restrict to the class of maps such that the image of 
the domain finite product group in each simple
component group of the codomain finite product
is acted on trivially by some automorphism arising from an
odd permutation. Two different maps of this kind, with the same domain
and codomain and the
same multiplicities, which can in any case differ at most by some 
automorphism leaving each component of the codomain invariant
and so arising there from a permutation, must then in fact
differ by an even permutation in each component---i.e., must
differ by an inner automorphism.
This happens for instance if one of the
components of the domain group is on an odd number of symbols and
its multiplicity is at least two, or also if there are at
least two symbols for the codomain component group which
are fixed by the image of the domain group.

There is also a way to apply Theorem 3 to inductive limits of
the finite product groups under consideration without
restricting to
special maps in the inductive limit construction. If we
enlarge the specifed normal subgroups
of the automorphism groups in the statement of
Theorem 3---somewhat misleadingly (in the present case) referred to as the inner
automorphism groups!---to include certain non-inner automorphisms
arising naturally in the present case, namely,
the product automorphisms considered above in the
case of a finite product of alternating groups (on 
five or more symbols),
arising from a permutation on each component,
and the natural extensions of these to the inductive limit group
(these automorphisms can be characterized without reference to
a particular inductive limit decomposition),
then the hypotheses of Theorem 3 are still satisfied. In particular,
the composition of an arbitrary homomorphism with such an automorphism on the domain side
is equal to the composition of this same
homomorphism with another such automorphism on the codomain side---the
basic algebraic property of inner automorphisms still obtains.
  
In short, in either setting (restricted sequences and inner
automorphisms, or arbitrary sequences and generalized inner
automorphisms as described above),
the maps between the finite product building blocks
modulo the special automorphisms considered are always determined by the
multiplicities of the partial maps.

In
other words the classification category for the 
category of building blocks, i.e., the finite products of
alternating groups on 
five or more symbols (slightly restricted,
in the case that Theorem 3 is applied with actual inner automorphisms)
is exactly
the same as that described in Section 2 for the
category of finite direct sums of matrix algebras over
the complex numbers (restricted to those of order
five or more for the present comparison). (Multiplicity zero
for a map between full matrix algebras means it is the
zero *-algebra map, while for a map between alternating groups
it means it is the trivial group map, and the analogy
is also close for higher multiplicities.)

It follows that the classifying category for
(sequential) inductive limits of the building block groups under
consideration
(finite products of alternating groups, on 
five or more symbols)
is the same as Bratteli's classifying
category for AF algbras, described briefly at the beginning
of Section 3 (Bratteli diagrams). (More precisely, we must
restrict consideration to AF algebras constructed using only
matrix algebras of order
five or more, but up to stable
isomorphism this is everything.) (Incidentally, it
might be interesting to consider whether there is an
analogue for groups, at least for those in the present
class, of stable isomorphism for C*-algebras.)

It is interesting to note that this category is equivalent
to the category of (countable) dimension groups, i.e., (countable)
unperforated
order groups with the Riesz decomposition property---with
a specified upward directed downward hereditary generating
subset of
the positive cone---sometimes called a scale. (See [8], [5] and
[13].)

Another point worthy of note is that, with this
identification of the classifying category (common for
AF algebras and the category of groups described above),
whereas the classification functor in the C*-algebra
setting is a familiar one---namely, K$_0$---the
corresponding functor in the group setting---mapping into
the same class of ordered groups---would not seem
hitherto  to have been considered. (And can it 
even be defined directly?)

(Added November 15, 2007: After this paper was submitted for publication the
author discovered the article
by Y.~Lavrenyuk and N.~Nekrashevych, On classification of inductive limits
of direct products of alternating groups, {\it J.~London Math.~Soc.} {\bf 75}
(2007), 146--162, which shows that the groups considered above are classified
up to isomorphism by their (equivalence classes of) Bratteli diagrams. 
To be more precise, these authors consider only the groups arising from a slightly
restricted class of sequences, which would seem to be restricted in a somewhat
different way from the first---restricted---class considered above. Note that
the second class of sequences considered
above is not restricted.
A second point of difference, also minor, is that, instead
of observing that mappings between alternating groups on five or more symbols
are automatically diagonal, with a certain multiplicity, up to conjugacy
by an inner automorphism, these authors consider only diagonal maps.)

\medskip

{\bf 5.2.}\ \
As pointed out in [8], the classification of 
AF algebras---the C*-algebra inductive limits of finite-dimensional
C*-algebras (i.e., finite C*-algebra direct sums of full matrix
algebras over the complex numbers)---is in a certain
sense equivalent to that of the corresponding *-algebra
inductive limits, or even just of the corresponding algebra
inductive limits. This sense can be
extended to the present context as follows: while these
three categories are presumably not equivalent, their
classifying categories given by Theorem 3 are all
equivalent---and are equivalent to the category just described
(which is sometimes referred to as the category of scaled
dimension groups).
\medskip

{\bf 5.3.}\ \
A relatively simple concrete
identification of the classifying category of Theorem 3,
for other subcategories of  the category of separable
C*-algebras considered
in Section 4.3 than the category of AF algebras
just discussed, would of course be interesting. To a
considerable extent, this is in fact
how the program of classifying (various classes of)
amenable C*-algebras has proceeded so far. (One
would perhaps like to consider the class of all amenable
C*-algebras at once---but not only 
did one case
take longer in the analogous setting of amenable
von Neumann algebras, some amenable
C*-algebras will definitely be more difficult than others---see
[33] and [32].)

For instance, if one considers simple inductive limits
of matrix algebras over C$([0, 1])$, say assumed to be
unital, then it follows from [10] that the classification
category of Theorem 3 consists of the category of order-unit
ordered groups arising in the (simple, unital) AF
case, with the modification that each one
should be paired at the same time with a Choquet
simplex (arising as the simplex of tracial states of the C*-algebra),
with the maps respecting this pairing.

In [16], the same
invariant---K$_0$ plus traces---, augmented by K$_1$, was
shown to be complete when the interval [0, 1] is
replaced by an arbitrary (variable) compact metric space
of  dimension 
at most three (or, in fact, any fixed number, but
with no new examples appearing). However, as was shown
already in [23], maps between C*-algebras
are not determined by these invariants up to
approximate unitary equivalence---even in the case of the
circle one needs to
consider, instead of the Banach algebra
K$_1$-groups, the (Hausdorffized) algebraic K$_1$-group 
(see [23] and [11]). In
the case of arbitrary compact metric spaces of 
dimension at most three, one must also
consider the K-groups with coefficients introduced
in [3] (and considered in the non-simple
real rank zero case in [6] and [2]). In fact, these
invariants suffice, as can be seen by study of [16] (and
can be seen immediately from Theorem 8.6 and Lemma 6.9 of [22]!). 
(The same result also holds in the more general case
considered recently by Niu in [25] (see also [17]), as
follows from Theorem 6.2.3 of [25] together with Lemma
6.9 of [22].)

In other words, for the class of C*-algebras classified
in [16] (or, more generally, in [25] and [17]), the abstract
classification functor of Theorem 3
is equivalent (by means of an equivalence of categories)
to the standard K-theoretical functor consisting of the  
invariants just listed. Furthermore, as shown in [34],
the objects arising as the values of the functor,
for the class of C*-algebras considered in [16], can be
described in simple terms (much as in the simpler cases
reported above). Incidentally, if the algebraic
K$_1$-group is taken to be based on invertible elements
rather than unitaries, then the invariant consisting
of the simplex of tracial states becomes redundant, as
the larger K$_1$-group is just the direct sum of the one
based on unitaries and the group of continuous affine
real-valued functions on the simplex.
(The affine
function corresponding to a projection is seen in this
picture as just the K$_1$-class of the exponential of the projection.)



There is another case in which the answer 
is simpler!---only the abstract K$_0$- and K$_1$-groups
and K-groups with coefficients. This is the case of
Kirchberg algebras (simple purely infinite separable
amenable C*-algebras) classified by Kirchberg and Phillips
in [19] and [26]. The classification functor of Theorem 3
is characterized as the functor KL of Rordam
([27]). (The range of this is still an abstract category, 
but in the case that the algebras satisfy
the Universal Coefficient Theorem (possibly automatic)---see
[30] and [4]---this is equivalent to the concrete
category of K-groups
with coefficients (including of course K$_0$ and K$_1$) 
referred to above.
\medskip

{\bf 5.4.}\ \
Consider the category of countably generated
Hilbert modules over a given C*-algebra A, with
embeddings, as in Section 4.4. (Recall that by Theorem 3.5
of [20], an $A$-module map between Hilbert
$A$-modules is an embedding in the present sense---i.e.,
preserves the $A$-valued inner product---if and only if
it is isometric.) In general, the structure
of the classifying category of Theorem 3---in particular,
the question when two morphisms are the same---would
appear to be somewhat complicated.
This can be seen already with the failure of cancellation for
isomorphism classes of algebraically finitely generated 
projective modules.

Remarkably, in the case that $A$ has stable
rank one (i.e., the ring $A$ with unit adjoined has Bass
stable rank one), the structure of the classifying category
becomes extremely simple: between any two objects, either
there are no maps, or there is exactly one map. (In other
words, any two morphisms between Hilbert C*-modules
are approximately equal modulo inner automorphisms; this
is proved in the last paragraph of the proof of Theorem 3 of [14].)
(This category is then
very much like the (common) classifying category given by 
Theorem 1 for sets, vector spaces,
or Hilbert spaces, with injective maps as morphisms 
(isometries in the
third case), and with the whole automorphism group
taken as the specified normal subgroup---namely,
cardinal numbers with maps just the relations $a\le b$.)

As a consequence (just as for sets!), a
Cantor-Bernstein theorem holds for
the category of countably generated Hilbert $A$-modules
(with embeddings) in the case that $A$ has stable rank one.
Indeed, if one has maps
$a\to b$ and $b\to a$ between objects $a$ and $b$, then these
persist in the classifying category, and
by uniqueness the composed maps in the classifying category
must be the identities 
for $a$ and $b$. In other words, the objects $a$ and $b$
are isomorphic in the classifying category, and hence by
Theorem 3---without using any more that
$A$ has stable rank one---they are isomorphic as Hilbert $A$-modules.
\medskip

{\bf 5.5.}\ \
Consider the category of pointed von Neumann
algebras with separable predual, as described in
Section 4.5. A certain subcategory of this is of particular
interest---namely, that for which the maps, in addition
to intertwining the specified states, and  their
modular automorphism
groups---and taking the centre into the centre (automatic
in the factor case)---so that
they pass to the Takesaki one-parameter crossed products---take 
the centre into
the centre at the level of the crossed products. The
subcategory comprising these maps admits a functor
to the Connes-Takesaki flow of weights---which is
just the restriction to the centre of the dual $\IR$-action
on the Takesaki crossed product. In
the amenable case this functor is, famously, a
classification functor. (To be precise one must restrict
to the case of properly infinite and continuous von Neumann
algebras---i.e., those of type II$_\infty$ or III. For these 
algebras the flow determines the algebra, and furthermore,
the flow may be an arbitrary measurable flow---by which is meant
an $\IR$-action on an
abelian von Neumann algebra with separable pre-dual.)

In the general case, the flow of weights functor
factors through the classification functor of Section 4.5---as
inner automorphisms belonging to the category of Section 4.5
also belong to the present subcategory, and their action on
the flow of weights is trivial---so that two morphisms in the
present subcategory which are approximately unitarily
equivalent, with respect to inner automorphisms in the
category under consideration---i.e., which are equal
in $\overline{\cC^{\out}}$---give
rise to the same morphism at the level of the
flows of weights.

Is this functor, from $\overline{\cC^{\out}}$
to the category of flows, in fact an equivalence of categories
in the amenable case? It would appear that this might follow from
careful inspection of [18], which shows that at least this is the case if one
considers only isomorphisms. (It remains to check that, also in the setting of
homomorphisms---between pointed von Neumann algebras---, two morphisms
between two objects in this category are approximately unitarily equivalent,
with respect to inner automorphisms in this category,
if they agree at the level of flows of weights.)
 
It should be noted that a general question that arises in the
setting of classification functors, and in the setting of
Theorem 3 in particular, is whether a functor exists in
the opposite direction, which is a one-sided inverse to the
classification functor. (One does not expect a two-sided
inverse, and in the setting of Theorem 3, at least, it
cannot exist, except in the trivial case that there
are no inner automorphisms---as the
classification functor kills all inner automorphisms.)
In [31], remarkably, a one-sided inverse is shown to exist
for the category of amenable pointed von Neumann algebras
with separable predual---in the restricted setting of
isomorphisms only---both at the level of the algebras and at
the level of the invariant, which as observed above is
the same whether it is considered as the abstract
category of Theorem 3, or as the category of flows.

The following concrete consequence of 
Theorem 3 in the case of the
category of Section 4.5 (so far, the only one!) is of
modest interest: namely,
restricting to finite factors (pointed with respect to
the unique tracial state) one obtains immediately the
Murray-von Neumann uniqueness theorem for the
approximately finite-dimensional case---or at least the 
important special case that there is a locally
finite-dimensional generating sub-*-algebra. (Any two maps
from such a finite factor
into an arbitrary finite factor are very easily seen
to be approximately unitarily equivalent---as the
calculation reduces immediately to the finite-dimensional
case for the domain algebra---and, furthermore,
at least one such map---constructed inductively---is
easily seen to exist. One is
therefore in the situation of a subcategory with a classifying
category, between any two objects of which there is exactly
one map!
(Necessarily, of course, an isomorphism.))
\bigskip

{\bf References}
\medskip
\itemitem{1.\ }   
O. Bratteli, {\it Inductive limits of finite dimensional C*-algebras,}
Trans. Amer. Math. Soc. {\bf 171} (1972), 195--234.
\smallskip

\itemitem{2.\ } 
M. Dadarlat and G. Gong, {\it A classification result for
approximately homogeneous C*-algebras of real rank zero}, Geom.
Funct. Anal. {\bf 7} (1997), 646--711.
\smallskip

\itemitem{3.\ } 
M. Dadarlat and T. A. Loring, {\it Classifying C*-algebras
via ordered mod-p $K$-theory}, Math. Ann. {\bf 305}
(1996), 601--616.
\smallskip

\itemitem{4.\ }
M. Dadarlat and T. A. Loring, {\it A universal multicoefficient
theorem for the Kasparov groups}, Duke Math. J. {\bf 24}
(1996), 355--377.
\smallskip

\itemitem{5.\ } 
E.~G. Effros, D.~E. Handelman, and C.-L. Shen,
{\it Dimension groups and their affine representations},
Amer. J. Math {\bf 102} (1980), 385--407.
\smallskip

\itemitem{6.\ }
S. Eilers, {\it A complete invariant for AD algebras with real
rank zero and bounded torsion in $K_1$,}
J. Funct. Anal. {\bf 139} (1996),  325--348.
\smallskip

\itemitem{7.\ }  
S. Eilers and G. Restorff, {\it On R\o rdam's classification
of certain C*-algebras with one non-trivial ideal,} Operator Algebras:
The Abel Symposium 2004 (editors, O. Bratteli, S. Neshveyev, and C.
Skau), Abel Symposia, {\bf 1}, Springer, Berlin, 2006, pages 87--96.
\smallskip

\itemitem{8.\ }
G. A. Elliott, {\it On the classification of inductive limits
of sequences of semisimple finite-dimensional algebras},
J. Algebra {\bf 38} (1976), 29--44.
\smallskip

\itemitem{9.\ }
G. A. Elliott, {\it On the classification of C*-algebras of real
rank zero}, J. Reine Angew. Math. {\bf 443} (1993), 179--219.
\smallskip

\itemitem{10.\ }
G.A. Elliott, {\it A classification of certain simple C*-algebras,}
Quantum and Non-Commutative Analysis 
(editors, H. Araki et al.),
Kluwer, Dordrecht, 1993, pages 373--385.
\smallskip

\itemitem{11.\ } 
G. A. Elliott, {\it A classification of certain simple C*-algebras,
II}, J. Ramanujan Math. Soc {\bf 12} (1997),  97--134.
\smallskip

\itemitem{12.\ }  
G. A. Elliott, {\it The classification problem for amenable C*-algebras,}
Proceedings of the International Congress of Mathematicians,
Z\"urich, 1994 (editor, S. D. Chatterji), Birkh\"auser, Basel, 1995, 
pages 922--932.
\smallskip

\itemitem{13.\ }
G. A. Elliott, {\it The inductive limit of a Bratteli diagram}, Lecture, MSRI,
Berkeley, September 2000. (Transparencies on MSRI web site.)
\smallskip

\itemitem{14.\ }
G. A. Elliott, K. T. Coward, and C. Ivanescu,
{\it The Cuntz semigroup as an invariant for C*-algebras,}
preprint.
\smallskip

\itemitem{15.\ } 
G. A. Elliott, D. E. Evans, and H. Su, {\it  Classification of inductive limits of matrix 
algebras over the Toeplitz algebra,} 
 Operator Algebras and Quantum Field Theory
(editors, S. Doplicher et al.), Accademia Nazionale dei Lincei, 
Rome, 1998, pages 36--50.
\smallskip

\itemitem{16.\ } 
G. A. Elliott, G. Gong, and L. Li, {\it On the classification of
simple inductive limit C*-algebras II: The isomorphism theorem,}
preprint, 1998 (to appear, Invent. Math.).
\smallskip

\itemitem{17.\ }
G. A. Elliott and Z. Niu, {\it On tracial approximation},
preprint.
\smallskip

\itemitem{18.\ }
Y. Kawahigashi, C. E. Sutherland, and M. Takesaki,
{\it The structive of the automorphism group of an injective
factor and the cocycle conjugacy of discrete abelian group
actions}, Acta Math. {\bf 169} (1992), 105--130.
\smallskip

\itemitem{19.\ }
E. Kirchberg, {\it The classification of purely infinite C*-algebras
using Kasporov's theory,} preprint, 1994.
\smallskip

\itemitem{20.\ }
E. C. Lance, {\it Hilbert C*-modules. A tool kit for operator 
algebraists.}  London Mathematical Society Lecture Note Series,
{\bf 210}, Cambridge University Press, Cambridge, 1995.
\smallskip

\itemitem{21.\ }
H. Lin, {\it An introduction to the classification of
amenable C*-algebras}, World Scientific Publishing Co., Inc.,
River Edge, NJ, 2001.
\smallskip

\itemitem{22.\ }
H. Lin, {\it Simple nuclear C*-algebras of tracial topological rank
one}, preprint.
\smallskip

\itemitem{23.\ }
K. E. Nielsen and K. Thomsen, {\it Limits of circle algebras,}
Expo. Math. {\bf 14} (1996), 17--56.
\smallskip

\itemitem{24.\ }
Z. Niu, {\it On the classfication of TAI algebras},
C. R. Math. Acad. Sci. Soc. R. Can. {\bf 26} (2004), 18--24.
\smallskip

\itemitem{25.\ }
Z. Niu, {\it On the classification of certain tracially
approximately subhomogeneous C*-algebras,} Ph.D. thesis,
University of Toronto, 2005.
\smallskip

\itemitem{26.\ }
N. C. Phillips, {\it A classification theorem for nuclear purely 
infinite simple C*-algebras}, Documenta Math. (2000), 
49--114.
\smallskip

\itemitem{27.\ }
M. R\o rdam, {\it Classification of certain infinite simple 
C*-algebras}, J. Funct. Anal. {\bf 131} (1995), 415--458.
\smallskip

\itemitem{28.\ }
M. R\o rdam, {\it Classification of extensions of certain 
C*-algebras by their six term exact sequences in $K$-theory,}
Math. Ann. {\bf 308} (1997), 97--117.
\smallskip

\itemitem{29.\ }
M. R\o rdam, {\it Classification of nuclear, simple
C*-algebras,} pages 1--145 of Classification of nuclear
C*-algebras. Entropy in operator algebras.
Encyclopedia of the Mathematical Sciences, {\bf 126},
Springer, Berlin, 2002.
\smallskip

\itemitem{30.\ }
J. Rosenberg and C. Schochet, {\it The K\"unneth Theorem and the
Universal Coefficient theorem for Kasparov's generalized $K$-functor,}
Duke Math. ~J. {\bf 55} (1987),  431--474.
\smallskip

\itemitem{31.\ }
C. Sutherland and M. Takesaki, {\it Right inverse
of the module of approximately finite dimensional
factors of type III and approximately finite
ergodic principal measured groupoids.}
Operator Algebras and Their Applications II
(editors, P. A. Fillmore and J. A. Mingo), 
Fields Institute Communications, {\bf 20},
American Mathematical Society, Providence, RI,
1998, pages 149--159.
\smallskip

\itemitem{32.\ }
A. Toms, {\it On the classification problem for
nuclear C*-algebras,} Ann. ~of Math., to appear.
\smallskip

\itemitem{33.\ }
J. Villadsen, {\it Simple C*-algebras with perforation},
J. Funct. Anal. {\bf 154} (1998),  110--116.
\smallskip

\itemitem{34.\ }
J. Villadsen, {\it The range of the Elliott invariant of the simple AH-algebras with
slow dimension growth}, K-Theory {\bf 15} (1998), 1--12. 

\vglue .5truein

Department of Mathematics,\par
University of Toronto,\par
Toronto, Canada\ \ M5S 2E4\par
\smallskip

\end

\bigskip

\vfill\eject

\end